\documentclass{amsart}

\usepackage{amsmath,amsfonts,amssymb,amscd}

\theoremstyle{definition}

\title{Division by $2$ on elliptic curves}
\author {Yuri G. Zarhin}

\address{Pennsylvania State University, Department of Mathematics, University Park, PA 16802, USA}

\email{zarhin@math.psu.edu}

\begin{document}\date{}

\maketitle
Let $K$ be a field of characteristic different from $2$.
Let 
$$E: y^2=(x-\alpha_1)(x-\alpha_2)(x-\alpha_3)$$
be an elliptic curve over $K$ where $\alpha_1,\alpha_2,\alpha_3$ are distinct elements of $K$.
Let $P=(x_2,y_2)$ be a $K$-point on $E$. It is well known that $P$ is divisible by 2 in $E(K)$ if and only if all three $x_2-\alpha_i$ are squares in $K$
\cite[p. 269--270]{Cassels}, 
\cite[Ch. 5, pp. 102--104]{Lang},  \cite{Husemoller}, \cite[Th. 4.2 on pp. 85-87]{Knapp}, \cite[pp. 331--332]{Bombieri}, \cite[pp. 212--214]{Wash}; see also \cite{Yelton}). 
This assertion is traditionally used in the course of a proof of the Weak Mordell-Weil Theorem for elliptic curves.
While the proof of the claim that the divisibility implies the squareness is straightforward, 
it seems that the known elementary  proofs of the converse statement are 
more involved/computational. (Notice that there is another approach that is based on Galois cohomology \cite[Sect. X.1, pp. 313--315]{Silverman}.)

Here we suggest an elementary proof of the divisibility that seems to be less computational. So, let us assume that all three $x_2-\alpha_i$ are squares in $K$.
Let $Q=(x_1,y_1)$ be a point on $E$ with $2 Q=P$.  Since $P\ne \infty$, $y_1\ne 0$ and therefore
the equation of the {\sl tangent line} $L$ to $E$ at $Q$ may be written in the form  
$$L: y=l x+m.$$
(Here $x_1,y_1, l, m$ are elements of an  overfield of $K$.) In particular,
$y_1=l x_1+m$.
By definition of $Q$ and $L$, the point $-P=(x_2,-y_2)$ is the ``third'' common point of $L$ and $E$;
in particular,
$-y_2=l x_2+m$, 
i.e., $y_2=-(l x_2+m)$.
Standard arguments (the restriction of the equation for $E$ to $L$, see \cite[pp. 25--27]{SilvermanTate}, \cite[pp. 12--14]{Wash} \cite[p. 331]{Bombieri}) tell us that
the cubic polynomial 
$$(x-\alpha_1)(x-\alpha_2)(x-\alpha_3)-(l x+m)^2$$
coincides with $(x-x_1)^2 (x-x_2)$. This implies that
$$-(l \alpha_i+m)^2=(\alpha_i-x_1)^2 (\alpha_i-x_2) \ \forall i=1,2,3.$$
Since $2Q=P\ne \infty$, none of   $x_1 -\alpha_i$ vanishes. Recall that all  $x_2-\alpha_i$ are squares in $K$ and they are obviously distinct. This implies that 
 corresponding square roots \cite[p. 331]{Bombieri}
$$r_i:=\frac{l\alpha_i+m}{x_1 - \alpha_i}=\sqrt{x_2-\alpha_i}$$
are {\sl distinct} elements of $K$. In other words, the transformation of the projective line 
$$z \mapsto \frac{l z+m}{-z+x_1}$$
sends three distinct $K$-points $\alpha_1,\alpha_2,\alpha_3$ to three distinct 
$K$-points $r_1,r_2,r_3$ respectively.
This implies that our transformation is {\sl not} constant (i.e., is an honest fractional-linear transformation) and defined over $K$.
Since one of the ``matrix entries",
$-1$ is already a nonzero element of $K$,  all other matrix entries $l, m, x_1$ also lie in $K$. Since $y_1=l x_1 +m$, it also lies in $K$.
So, $Q=(x_1,y_1)$ is a $K$-point of $E$.

Let us get explicit formulas for $x_1,y_1, l, m$ in terms of $r_1,r_2,r_3$. We have
$$\alpha_i=x_2-r_i^2, \ l\alpha_i+m=r_i (x_1 - \alpha_i)$$
and therefore
$$l (x_2-r_i^2)+m=r_i [x_1-(x_2-r_i^2)]=r_i^3+(x_1-x_2)r_i,$$
which is equivalent to
$r_i^3+ l r_i^2+(x_1-x_2)r_i-(l x_2+m)=0$
and this equality holds for all $i=1,2,3$. This means that the (monic) cubic polynomial
$$h(t)=t^3+l t^2+(x_1-x_2)t-(l x_2+m)$$
coinsides with $(t-r_1)(t-r_2)(t-r_3)$. Recall that $-(l x_2+m)=y_2$ and get
$$r_1 r_2 r_3=-y_2.$$
We also get
$$l=-(r_1+r_2+r_3), \ x_1-x_2=r_1 r_2+r_2 r_3+r_3 r_1,$$
This implies that
$$x_1=x_2+(r_1 r_2+r_2 r_3+r_3 r_1).$$
Since $y_1=l x_1+m$ and $-y_2=l x_2+m$, we obtain that
$$m=-y_2-l x_2=-y_2+(r_1+r_2+r_3)x_2,$$
and therefore
$$y_1= -(r_1+r_2+r_3)[x_2+(r_1 r_2+r_2 r_3+r_3 r_1)]+[-y_2+(r_1+r_2+r_3)x_2],$$ i.e.,
$$y_1=-y_2-(r_1+r_2+r_3)(r_1 r_2+r_2 r_3+r_3 r_1).$$
Notice that there are precisely four points $Q \in E(K)$ with $2Q=P$, each of which corresponds to one of {\sl four} choices of three square roots $r_i=\sqrt{x_2-\alpha_i}\in K$ ($i=1,2,3$) with
$r_1 r_2 r_3=-y_2$
in such a way  that the corresponding 
$$Q=\left(x_2+(r_1 r_2+r_2 r_3+r_3 r_1),-y_2-(r_1+r_2+r_3)(r_1 r_2+r_2 r_3+r_3 r_1)\right).$$

{\bf Example}.
 Let us choose as $P=(x_2,y_2)$ the point $(\alpha_3,0)$  of order $2$ on $E$, Then $r_3=0$ and we have two arbitrary independent choices of  (nonzero) $r_1=\sqrt{\alpha_3-\alpha_1}$ and $r_2=\sqrt{\alpha_3-\alpha_2}$.
 Then 
 $$Q=(\alpha_3+r_1 r_2, -(r_1+r_2)r_1 r_2)=(\alpha_3+r_1 r_2, -r_1(\alpha_3-\alpha_2)-r_2(\alpha_3-\alpha_1))$$
 is a point on $E$ with $2Q=P$; in particular, $Q$ is a point of order $4$. The same is true for (three remaining) points
 $-Q=(\alpha_3+r_1 r_2, r_1(\alpha_3-\alpha_2)+r_2(\alpha_3-\alpha_1))$, 
 \newline
 $(\alpha_3-r_1 r_2, - r_1(\alpha_3-\alpha_2)+r_2(\alpha_3-\alpha_1))$ and 
 $(\alpha_3-r_1 r_2,  r_1(\alpha_3-\alpha_2)-r_2(\alpha_3-\alpha_1))$.


{\bf Acknowledgements}. This work was partially supported by a grant from the Simons Foundation (\#246625 to Yuri Zarkhin). This note was written during my stay in May-June 2016 at the Max-Planck-Institut f\"ur Mathematik (Bonn, Germany),
whose hospitality and support are gratefully acknowledged.

\end{document}